\documentclass[10pt]{article}

\usepackage[utf8]{inputenc}
\usepackage{amsmath}
\usepackage{amsfonts}
\usepackage{amssymb}
\usepackage{eufrak}
\usepackage{mathrsfs}
\usepackage{dsfont}
\usepackage{fancyhdr}
\usepackage{indentfirst}
\usepackage{graphicx}
\usepackage{newlfont}
\usepackage{amssymb}
\usepackage{amsmath}
\usepackage{latexsym}
\usepackage{amsthm}
\usepackage{mathtools}
\usepackage{nicefrac}
\usepackage{epstopdf}
\usepackage{caption}
\usepackage{color}
\usepackage{dsfont}
\usepackage{comment}
\usepackage{hyperref}

\usepackage{esvect}

\usepackage{pdfsync}
\usepackage{amsmath}                                                        
\numberwithin{equation}{section}
\usepackage{graphicx}                                                       
\usepackage{subfigure}
\usepackage{enumitem}                                                    
\usepackage{hyperref}
\usepackage{verbatim}
\usepackage{xcolor}

\usepackage{amssymb}                                                        
\usepackage[mathscr]{eucal}                                             
\usepackage{cancel}                                                             
\usepackage[normalem]{ulem}                                                                 
\usepackage{pstricks}
\usepackage{rotating}
\usepackage{lscape}
\usepackage[paperwidth=8.5in,paperheight=11in,top=1.00in, bottom=1.00in, left=1.00in, right=1.00in]{geometry}
\usepackage{mathtools}                                                      
\mathtoolsset{showonlyrefs=true}                                    
\usepackage{fixltx2e,amsmath}                                           
\MakeRobust{\eqref}

\usepackage{dsfont}

\usepackage{mathdots}
\usepackage{amsthm}                                                             
\allowdisplaybreaks

\theoremstyle{plain}
\newtheorem{theorem}{Theorem}
\numberwithin{theorem}{section}

\newtheorem{lemma}[theorem]{Lemma}                              
\newtheorem{proposition}[theorem]{Proposition}
\newtheorem{corollary}[theorem]{Corollary}
\theoremstyle{definition}

\newtheorem{remark}[theorem]{Remark}
\newtheorem{assumption}[theorem]{Assumption}

\makeatletter
\DeclareRobustCommand{\cev}[1]{%
  {\mathpalette\do@cev{#1}}%
}
\newcommand{\do@cev}[2]{%
  \vbox{\offinterlineskip
    \sbox\z@{$\m@th#1 x$}%
    \ialign{##\cr
      \hidewidth\reflectbox{$\m@th#1\vec{}\mkern4mu$}\hidewidth\cr
      \noalign{\kern-\ht\z@}
      $\m@th#1#2$\cr
    }%
  }%
} \makeatother

\def \s {{\sigma}}

\def \a {{\alpha}}
\def \b {{\beta}}

\def \R {\mathbb{R}}
\def \p {\partial}

\def \PP {\mathsf{P}}

\newcommand\CC{\mathcal{C}}

\newcommand\dd{\mathrm{d}}

\def \R  {{\mathbb {R}}}

\def \n {{\nu}}
\def \m {{\mu}}

\def \z {{\zeta}}
\def \p {{\partial}}
\def \a {{\alpha}}
\def \O {{\Omega}}

\def \a {{\alpha}}
\def \b {{\beta}}

\def \mm {\mathfrak{m}}
\def \MM {\mathfrak{M}}

\def \It\^o {It\^o }
\def \s {{\sigma}}

\def \R {{\mathbb {R}}}

\def \n {{\nu}}
\def \v {{\nu}}
\def \m {{\mu}}

\def \O {{\Omega}}
\def \phi {{\varphi}}

\def\l {\lambda}

\def \A {\mathcal{A}}

\def \F {\mathcal{F}}

\def \Ã  {{\`a }}
\def \Ã¨ {{\`e }}
\def \Ã² {{\`o }}
\def \Ã¹ {{\`u }}

\begin{document}

\title{McKean-Vlasov stochastic equations with H\"older coefficients}

\author{
Andrea Pascucci\thanks{Dipartimento di Matematica, Universit\`a di Bologna, Bologna, Italy.
\textbf{e-mail}: andrea.pascucci@unibo.it} \and Alessio Rondelli\thanks{Dipartimento di
Matematica, Universit\`a di Bologna, Bologna, Italy. \textbf{e-mail}:
alessio.rondelli@studio.unibo.it}}

\date{This version: \today}

\maketitle

\begin{abstract}
This work revisits the well-posedness of non-degenerate McKean-Vlasov stochastic differential
equations with H\"older continuous coefficients, recently established by Chaudru de Raynal. We
provide a streamlined and direct proof that leverages standard Gaussian estimates for uniformly
parabolic PDEs, bypassing the need for derivatives with respect to the measure argument and
extending applicability to hypoelliptic PDEs under weaker assumptions.
\end{abstract}

\section{Introduction}\label{intro}
In \cite{MR4035024}, Chaudru~de Raynal recently established the well-posedness of non-degenerate
McKean-Vlasov (MKV) stochastic differential equations with H\"older drift. His approach is based
on the associated PDE formulated on the domain $[0,T] \times \mathbb{R}^d \times
\mathcal{P}^2(\mathbb{R}^d)$, where $T
> 0$, $d$ is the dimensionality of the system, and $\mathcal{P}^2(\mathbb{R}^d)$ represents the
space of probability measures on $\mathbb{R}^d$ with finite second moment.

This note aims to present a more streamlined and direct proof of a more general result, avoiding
the use of PDEs incorporating derivatives with respect to the measure variable. Our proof relies
solely on standard Gaussian estimates for uniformly parabolic PDEs on $[0,T] \times \mathbb{R}^d$.
This approach not only simplifies the argument but also renders the proof adaptable to broader
classes of degenerate MKV equations, where the associated PDEs are hypoelliptic, provided that
upper Gaussian bounds for their fundamental solutions are available.

Let $\mathcal{P}(\mathbb{R}^d)$ denote the space of probability measures on $\R^{d}$ and $[Z]$ be
the law of a random variable $Z$. We consider the MKV equation
\begin{equation}\label{e1}
  X_{t}=X_{0}+\int_{0}^{t}B(s,X_{s},[X_{s}])\dd s+\int_{0}^{t}\Sigma(s,X_{s},[X_{s}])\dd W_{s},\qquad
  t\in[0,T],
\end{equation}
where $W$ is a standard $d$-dimensional Brownian motion defined on a
filtered probability space $(\O,\F,\PP,\F_{t})$ and $T>0$. 
We consider coefficients
\begin{equation}\label{e14}
  B:[0,T]\times\R^{d}\times\mathcal{P}(\mathbb{R}^d)\longrightarrow\R^{d},\qquad
  \Sigma:[0,T]\times\R^{d}\times\mathcal{P}(\mathbb{R}^d)\longrightarrow\R^{d\times d},
\end{equation}
of the form 
\begin{equation}\label{e2}
 B(t,x,\m)=\int_{\R^{d}}b(t,x,y)\m(\dd y),\qquad \Sigma(t,x,\m)=\int_{\R^{d}}\s(t,x,y)\m(\dd y).
\end{equation}
This structural assumption can be significantly weakened to allow for broader forms. For further
details and specific examples, refer to Remarks \ref{r4} and \ref{r5}.
\begin{assumption}[\bf
Regularity]\label{H1} The coefficients $b,\s\in L^{\infty}([0,T];bC^{\a}(\R^{d}\times\R^{d}))$ for
some $\a\in\,]0,1]$, where $bC^{\a}$ denotes the space of bounded and $\a$-H\"older continuous
functions, equipped with the norm
    $$\|f\|_{bC^{\a}}:=\sup_{z}|f(z)|+[f]_{C^{\a}},\qquad
    [f]_{C^{\a}}:=\sup_{z\neq\z}\frac{|f(z)-f(\z)|}{|z-\z|^{\a}}.$$
\end{assumption}
\begin{assumption}[\bf Non-degeneracy]\label{H2}
The matrix $\CC:=\Sigma\Sigma^{\ast}$ is uniformly positive definite, that is
\begin{equation}\label{e20}
  \langle \CC(t,x,\m)y,y\rangle\ge\l|y|^{2},\qquad t\in[0,T],\ x,y\in\R^{d},\
  \m\in\mathcal{P}(\mathbb{R}^d),
\end{equation}
for some positive constant $\l$.
\end{assumption}
Our main result is the following
\begin{theorem}[\bf Weak well-posedness]\label{t1}
Under Assumptions \ref{H1} and \ref{H2}, for any $\bar{\m}_{0}\in\mathcal{P}(\mathbb{R}^d)$ there
exists a unique weak solution of  \eqref{e1} such that $[X_{0}]=\bar{\m}_{0}$.
\end{theorem}
A direct consequence of Theorem \ref{t1} is the following
\begin{corollary}[\bf Strong well-posedness]\label{c1}
Under Assumptions \ref{H1} and \ref{H2}, if $\s=\s(t,x,y)$ is also Lipschitz continuous in $x$
uniformly in $(t,y)\in[0,T]\times\R^{d}$ then \eqref{e1} admits a unique strong solution.
\end{corollary}
\begin{remark}\label{r4}
Theorem \ref{t1} remains valid even if the coefficients $B,\Sigma$ in \eqref{e14} are not
necessarily of the form \eqref{e2}, provided they are bounded, satisfy Assumption \ref{H2} and
the H\"older condition
\begin{equation}\label{e16}
  |B(t,x,\m)-B(t,y,\n)|+|\Sigma(t,x,\m)-\Sigma(t,y,\n)|\le c
  (|x-y|^{\a}+\mm_{\a}(\m,\v)).
\end{equation}
Here $\mm_{\a}$ is the metric on $\mathcal{P}(\R^{d})$ defined by
\begin{equation}\label{e11}
 \mm_{\a}(\mu,\nu):=\sup_{\|f\|_{bC^\alpha}\leq 1}\int_{\R^{d}} f(x)\left(\mu-\nu\right)(\dd x).
\end{equation}
As \eqref{e16} ensures the validity of estimate \eqref{e13} in Section \ref{proof}, the proof of
Theorem \ref{t1} follows through without requiring significant modifications.
\end{remark}
\begin{remark}\label{r5}
Our approach can also easily handle the case where the metric $\mm_{\a}$ in 
\eqref{e16} is replaced by
 $$m_\a(\mu,\nu)=\sup_{[f]_{C^{\a}}\leq1}\int_{\R^{d}} f(x)(\mu-\nu)(\dd x),$$
for some $\a\in\,]0,1]$, provided that the initial distribution belongs to
$\mathcal{P}^\a(\R^{d})$. The result in \cite{MR4035024} fits as a special case. Specifically, in
\cite{MR4035024} the drift coefficient has the form
  $$B(t,x,\m)=\b\left(t,x,\int_{\R^{d}}\phi(y)\m(\dd y)\right)$$
where $\b=\b(t,x,w)$ is defined on $[0,T]\times \R^d\times \R$, is bounded, H\"older continuous in
$x$ uniformly in $(t,w)$, and is differentiable in $w$ with bounded derivative, uniformly in
$(t,x)$; additionally, $\phi\in C^{\a}(\R^{d})$ and the initial distribution is required to belong
to $\mathcal{P}^2(\R^d)$.
\end{remark}
As mentioned earlier, our approach is naturally suited for generalization to broader classes of
equations. Recently, in collaboration with A.Y. Veretennikov \cite{PRV24}, we studied the weak
well-posedness for class of {\it degenerate} MKV equations whose prototype is the kinetic system
\begin{equation}\label{eL}
  \begin{cases}
    dX_{t}=V_{t}dt, \\
    dV_{t}=B(t,X_{t},V_{t},[(X_{t},V_{t})])dt
    +\Sigma(t,X_{t},V_{t},[(X_{t},V_{t})])dW_{t}.
  \end{cases}
\end{equation}
The classical Langevin model \cite{Langevin} serves as a particular instance of \eqref{eL}, where
the solution is a $2d$-dimensional process $(X_{t}, V_{t})$ representing the dynamics of a system
of $d$ particles in phase space, with $X_{t}$ and $V_{t}$ denoting the position and velocity at
time $t$, respectively. Despite the failure of the non-degeneracy Assumption \ref{H2} in this
context and the associated Kolmogorov PDEs being hypoelliptic rather than uniformly parabolic, the
existence and Gaussian estimates for the fundamental solution were established by the first author
and collaborators in \cite{DiFrancescoPascucci2} and \cite{MR4660246}. This case seems attainable
and will be the subject of future work, where we aim to generalize the results of this paper to
degenerate MKV equations.

The remainder of the paper is structured as follows. Section \ref{proof2} introduces the
analytical tools necessary to frame a fixed-point problem in the space of continuous flows of
marginals. In particular, we prove a crucial lemma for the inversion of forward and backward
transport operators. Section \ref{proof} is dedicated to the proof of Theorem \ref{t1}. Finally,
the Appendix compiles some result on the completeness of spaces of measures utilized throughout
the paper.

\section{Inversion lemma}\label{proof2}
The following notations for integrals will be used interchangeably throughout this section:
  $$\int_{\R^{d}}
f\dd\n=\int_{\R^{d}} f(x)\n(\dd x)=\int_{\R^{d}} \n(\dd x)f(x).$$ For a fixed flow of marginals
$\m=(\m_{t})_{t\in[0,T]}\in C([0,T];\mathcal{P}(\mathbb{R}^d))$ and an initial distribution
$\bar{\m}_{0}\in \mathcal{P}(\mathbb{R}^d)$, we consider the ``linearized'' version of \eqref{e1}
\begin{equation}\label{e1b}
  \dd X_{t}=B(t,X_{t},\m_{t})\dd t+\Sigma(t,X_{t},\m_{t})\dd W_{t},\qquad
  [X_{0}]=\bar{\m}_{0},
\end{equation}
which is a standard (i.e. non-MKV) SDE and admits a unique weak solution denoted by $X^{\m}$.

The solution $X^{\m}$ is a Markov process with infinitesimal generator
\begin{equation}\label{e3}
  \A^{\m}_{t,x}:=\frac{1}{2}\sum_{i,j=1}^{d}\CC_{ij}^{\m}(t,x)\p_{x_{i}x_{j}}+\sum_{i=1}^{d}B_{i}^{\m}(t,x)\p_{x_{i}}
\end{equation}
where
\begin{equation}\label{e4}
  \CC_{ij}^{\m}(t,x):=(\Sigma\Sigma^{\ast})_{ij}(t,x,\m_{t}),\quad B_{i}^{\m}(t,x):=B_{i}(t,x,\m_{t}),\qquad (t,x)\in[0,T]\times\R^{d},\ 1\le i,j\le d.
\end{equation}
By the classical theory of uniformly parabolic PDEs with bounded and H\"older continuous
coefficients (see \cite{friedman-parabolic} or the more recent presentation in \cite{PascucciSC}),
operator $\p_{t}+\A^{\m}_{t,x}$ has a fundamental solution $p^{\m}=p^{\m}(t,x;s,y)$, defined for
$0\le t<s\le T$ and $x,y\in\R^{d}$, that is the transition density of $X^{\m}$.

The {\it push-forward} and the {\it pull-back} operators acting on the distribution
$\bar{\m}_{0}\in\mathcal{P}(\R^{d})$ are defined as
 $$
 \vec{P}^{\m}_{t,s}\bar{\m}_{0}(y)
 :=\int_{\R^{d}}p^\mu(t,x;s,y)\bar{\m}_{0}(\dd x),\qquad
 \cev{P}^{\m}_{t,s}\bar{\m}_{0}(x)
 :=\int_{\R^{d}} p^\mu(t,x;s,y)\bar{\m}_{0}(\dd y),$$
for $0\le t<s\le T$ and $x,y\in\R^{d}$. Notice that $\vec{P}^{\m}_{0,s}\bar{\m}_{0}$ is the
density of the marginal law $[X^{\m}_{s}]$.

\medskip 
A key element in the proof of Theorem \ref{t1} is the following inversion lemma, which expresses
the push-forward operator in terms of pull-back operators. This result draws inspiration from
\cite{MR3072241}, where a similar formula, equation (28) in \cite{MR3072241}, is presented, though
formulated in the ``opposite direction''.
\begin{lemma}[\bf Inversion]\label{l1}
Let $X^{\m},X^{\n}$ be the weak solutions of the linearized SDE \eqref{e1b} corresponding to the
flows of marginals $\m,\n\in C([0,T];\mathcal{P}(\mathbb{R}^d))$ respectively, and with the same
initial distribution $[X^{\m}_{0}]=[X^{\n}_{0}]=\bar{\m}_{0}\in\mathcal{P}(\R^{d})$. Then, for any
$f\in bC^{\a}(\R^{d})$ and $s\in\,]0,T]$, we have
\begin{equation}\label{e100}
  I_{s}^{\m,\n}(f):=\int_{\R^{d}}\dd y\,f(y)(\vec{P}^{\m}_{0,s}-\vec{P}^{\n}_{0,s})\bar{\m}_{0}(y)=
  \int_{\R^{d}}\bar{\m}_{0}(\dd x)\int_{0}^{s}\dd t\,\cev{P}^{\m}_{0,t}\left(\A^{\m}_{t,x}-\A^{\n}_{t,x}\right)\cev{P}^{\n}_{t,s}f(x).
\end{equation}
\end{lemma}
\proof Under Assumptions \ref{H1}, \ref{H2} and if the coefficients are also continuous with
respect to the time variable\footnote{This assumption simplifies the presentation without imposing
significant restrictions: even if the coefficients are merely measurable in $t$,
$p^{\m}(\cdot,x;s,y)$ remains absolutely continuous and the proof proceeds in a similar manner.},
it is well-known that the transition density $p^{\m}(\cdot,\cdot;s,y)\in
C^{1,2}([0,s[\times\R^{d})$ for any $(s,y)\in\,]0,T]\times\R^{d}$, and is a classical solution of
the backward Kolmogorov PDE
\begin{equation}\label{BKPDE}
  (\p_{t}+\A^{\m}_{t,x})p^{\m}(t,x;s,y)=0,\qquad (t,x)\in[0,s[\times\R^{d};
\end{equation}
moreover, for any $(t,x)\in[0,T[\times\R^{d}$, $p^{\m}(t,x;\cdot,\cdot)$ is a distributional
solution of the forward Kolmogorov PDE $(\p_{s}-\A^{\m,\ast}_{s,y})p^{\m}(t,x;s,y)=0$, that is
\begin{equation}\label{BKPDE}
  \int_{0}^{T}\dd s\int_{\R^{d}}\dd y\,p^{\m}(t,x;s,y)(\p_{s}+\A^{\m}_{s,y})\phi(s,y)=0,
\end{equation}
for any test function $\phi\in C_{0}^{\infty}\in\,]t,T]\times\R^{d}$.

We have
\begin{align}
  I_{s}^{\m,\n}(f)&=
    \int_{\R^{d}}\dd y\,f(y)\int_{\R^{d}}\bar{\m}_{0}(\dd x)(p^\mu(0,x;s,y)-p^{\n}(0,x;s,y))=
\intertext{(by Fubini's theorem and standard Gaussian estimates for $p^{\m}$ and $p^{\n}$; see
Chapter 1, Section 6 in \cite{friedman-parabolic} or Theorem 20.2.5 in \cite{PascucciSC})}
  &=\int_{\R^{d}}\bar{\m}_{0}(\dd x)\int_{\R^{d}}\dd y\,f(y)(p^\mu(0,x;s,y)-p^{\n}(0,x;s,y))\\
  &=\int_{\R^{d}}\bar{\m}_{0}(\dd x)\int_{\R^{d}}\dd y\,f(y)\int_{0}^{s}\dd
  t\frac{d}{dt}\int_{\R^{d}}\dd z\,p^\mu(0,x;t,z)p^{\n}(t,z;s,y)=
\intertext{(by classical potential estimates; see Chapter 1, Section 3 in
\cite{friedman-parabolic} or Proposition 20.3.9 in \cite{PascucciSC})}
  &=\int_{\R^{d}}\bar{\m}_{0}(\dd x)\int_{\R^{d}}\dd y\,f(y)\int_{0}^{s}\dd t
  \int_{\R^{d}}\dd z\left((\p_{t}p^\mu(0,x;t,z))p^{\n}(t,z;s,y)+p^\mu(0,x;t,z)\p_{t}p^{\n}(t,z;s,y)\right)=
\intertext{(by the forward and backward Kolmogorov equations)}
  &=\int_{\R^{d}}\bar{\m}_{0}(\dd x)\int_{\R^{d}}\dd y\,f(y)\int_{0}^{s}\dd t
  \int_{\R^{d}}\dd z\,p^\mu(0,x;t,z)\left(\A^{\m}_{t,z}-\A^{\n}_{t,z}\right)p^{\n}(t,z;s,y)=
\intertext{(using the potential estimates once again)}
  &=\int_{\R^{d}}\bar{\m}_{0}(\dd x)\int_{0}^{s}\dd t
  \int_{\R^{d}}\dd z\,p^\mu(0,x;t,z)\left(\A^{\m}_{t,z}-\A^{\n}_{t,z}\right)\int_{\R^{d}}\dd y\,f(y)p^{\n}(t,z;s,y)\\
  &=\int_{\R^{d}}\bar{\m}_{0}(\dd x)\int_{0}^{s}\dd
  t\,\cev{P}^{\m}_{0,t}\left(\A^{\m}_{t,x}-\A^{\n}_{t,x}\right)\cev{P}^{\n}_{t,s}f(x).
\end{align}
\endproof

\section{Proof of Theorem \ref{t1}}\label{proof}
The proof of Theorem \ref{t1} relies on a contraction mapping principle in the space
$C([0,T];\mathcal{P}(\mathbb{R}^d))$ of continuous flows of distributions $(\m_{t})_{t\in[0,T]}$,
equipped with the metric
\begin{equation}\label{e12}
   \MM_{\a}(\m,\v):=\max\limits_{t\in[0,T]}\mm_{\a}(\m_{t},\n_{t}),
\end{equation}
where $\mm_{\a}$ is defined in \eqref{e11}. Specifically,
denoting $X^{\m}$ the solution of \eqref{e1b}, we claim that the map 
\begin{align}
  \m\longmapsto[X^{\m}_{t}]_{t\in[0,T]}
\end{align}
is a contraction on $C([0,T];\mathcal{P}(\mathbb{R}^d))$, at least for $T>0$ suitably small. The
thesis will readily follow from this assertion. Indeed, by Proposition \ref{l2},
$(C([0,T];\mathcal{P}(\mathbb{R}^d)),\MM_{\a})$ is a complete metric space and
therefore there would exist a unique fixed point $\bar{\m}\in C([0,T];\mathcal{P}(\mathbb{R}^d))$, such that $\bar{\m}=[X^{\bar{\m}}]$. 
Thus, $X^{\bar{\m}}$ is the unique weak (or strong, under the assumption of Corollary \ref{c1})
solution of the MKV equation \eqref{e1}.

Recalling the notation \eqref{e4} for the coefficients of the infinitesimal generators, we have
the preliminary estimate
\begin{align}\label{e5}
 \|\CC^{\m}-\CC^{\n}\|_{L^{\infty}([0,T]\times\R^{d})}+\|B^{\m}-B^{\n}\|_{L^{\infty}([0,T]\times\R^{d})}\le
 c\, \MM_{\a}(\mu,\nu)
\end{align}
for some positive constant $c$ which depends only on the norms of $b$ and $\s$ in
$L^{\infty}([0,T];bC^{\a}(\R^{d}\times\R^{d}))$. 
Estimate \eqref{e5} follows from
\begin{align}\label{e13}
  |B^{\m}(t,x)-B^{\n}(t,x)|=\left|\int_{\R^{d}}b(t,x,y)\left(\m_{t}-\n_{t}\right)(\dd
  y)\right|\le \|b(t,x,\cdot)\|_{bC^{\a}(\R^{d})}\mm_{\a}(\m_{t},\n_{t})
\end{align}
and the analogous estimate for $\CC$.

Now, we note that
\begin{equation}\label{e22}
  \MM_{\a}(\mu,\nu)=\max_{s\in[0,T]}\sup_{\|f\|_{bC^\alpha}\le 1}I_{s}^{\m,\n}(f)
\end{equation}
with $I_{s}^{\m,\n}(f)$ as in \eqref{e100}, and we have
\begin{align}
  |I_{s}^{\m,\n}(f)|&\le
  \int_{0}^{s}\|\cev{P}^{\m}_{0,t}(\A^{\m}_{t,\cdot}-\A^{\n}_{t,\cdot})\cev{P}^{\n}_{t,s}f\|_{L^{\infty}(\R^{d})}\dd
  t\\
  &\le\int_{0}^{s}\|(\A^{\m}_{t,\cdot}-\A^{\n}_{t,\cdot})\cev{P}^{\n}_{t,s}f\|_{L^{\infty}(\R^{d})}\dd t\le
\intertext{(by \eqref{e5})}
  &\le c\, \MM_{\a}(\mu,\nu)\int_{0}^{s}\max_{1\le i,j\le d}
  \left(\|\p_{x_{i}x_{j}}\cev{P}^{\n}_{t,s}f\|_{L^{\infty}(\R^{d})}+
  \|\p_{x_{i}}\cev{P}^{\n}_{t,s}f\|_{L^{\infty}(\R^{d})}\right)\dd t\le
\intertext{(by the Gaussian estimates for $p^{\n}$, adjusting the constant $c$ as
needed)}\label{e7}
  &\le c\, \MM_{\a}(\mu,\nu)\int_{0}^{s}\frac{1}{(s-t)^{1-\frac{\a}{2}}}\dd t\le
  c\, T^{\frac{\a}{2}} \MM_{\a}(\mu,\nu).
\end{align}
In conclusion, combining \eqref{e22} and \eqref{e7}, we obtain
\begin{equation}
  \MM_{\a}(\mu,\nu)
  \le c\, T^{\frac{\a}{2}} \MM_{\a}(\mu,\nu)
\end{equation}
which proves the thesis.

\section{Appendix}
We provide the proof of routine results for which we could not find an appropriate reference. We
recall the notations
\begin{align}
  \MM_{\a}(\m,\v)&=\max\limits_{t\in[0,T]}\mm_{\a}(\m_{t},\n_{t}),\qquad
  \mm_{\a}(\m_{t},\n_{t})=\sup_{\|f\|_{bC^\alpha}\leq 1}\int_{\R^{d}} f(x)\left(\mu_{t}-\nu_{t}\right)(\dd
  x),\\
  M_{\a}(\m,\v)&=\max\limits_{t\in[0,T]}m_{\a}(\m_{t},\n_{t}),\qquad m_{\a}(\m_{t},\n_{t})=
  \sup_{[f]_{C^{\a}}\leq1}\int_{\R^{d}} f(x)(\mu_{t}-\nu_{t})(\dd x),
\end{align}
\begin{proposition}\label{l2}
For any $\a\in\,]0,1]$, $\left(C([0,T];\mathcal{P}(\R^d)),\MM_{\a}\right)$ and
$\left(C([0,T];\mathcal{P}^{\a}(\mathbb{R}^d)),M_{\a}\right)$ are complete metric spaces.
\end{proposition}
\proof It suffices to show that $(\mathcal{P}(\R^d),\mm_\a)$ and $(\mathcal{P}^\a(\R^d),m_\a)$ are
complete metric spaces. We use the fact that $(\R^d,d_\a)$, where $d_\a(x,y):=|x-y|^\a$, is a
Polish space. Then, we notice that $\mm_\a$ is the bounded Lipschitz distance with respect to
$d_\a$, that is
  $$\mm_\a(\mu,\nu)=\sup_{\|f\|_{b\text{Lip}(d_\a)}\leq 1}\int_{\R^{d}} f(x)(\mu-\nu)(\dd x)=:\|\mu-\nu\|^{\ast}_{\text{BL}(d_\a)},
  \qquad \m,\n\in\mathcal{P}(\mathbb{R}^d).$$
By Theorem 8.10.43 in \cite{bogachev2007measure},
$(\mathcal{P}(\R^d),\|\cdot\|^{\ast}_{\text{BL}(d_\a)})=(\mathcal{P}(\R^d),\mm_\a)$ is complete.
This proves the first part of the statement.

Analogously, $m_\a$ is the $1$-Wasserstein metric on $\mathcal{P}^{1}(\R^d)$, where $\R^{d}$ is
equipped with $d_\a$: indeed, we have
  $$m_\a(\mu,\nu)
  =\sup_{[f]_{\text{Lip}(d_\a)}\leq 1}\int_{\R^{d}}
  f(x)(\mu-\nu)(\dd x)=:W_{1,d_\a}(\mu,\nu), \qquad \m,\n\in\mathcal{P}^{1}(\mathbb{R}^d,d_{\a}).$$
By Theorem 6.18 in \cite{MR2459454},
$(\mathcal{P}^1(\R^d,d_{\a}),W_{1,d_\a})=(\mathcal{P}^\a(\R^d),m_\a)$ is a complete metric space.
This concludes the proof.
\endproof

\bibliographystyle{plain}
\bibliography{Plain}

%
\end{document}